\documentclass[]{article}

\usepackage{lineno,hyperref,amsmath,amssymb,mathrsfs,color,graphicx}
\modulolinenumbers[5]

\newcommand{\bfx}{\mathbf{x}}
\newcommand{\RE}{\mathbb{R}}

\newtheorem{theorem}{Theorem}
\newtheorem{proposition}{Proposition}
\newtheorem{remark}{Remark}

\begin{document}


\title{A reduced order model for the finite element approximation of
eigenvalue problems}

\author{}
\author{}
\author{Fleurianne Bertrand, Daniele Boffi, Abdul Halim}
\date{}
\maketitle



\begin{abstract}
In this paper we consider a reduced order method for the approximation of the
eigensolutions of the Laplace problem with Dirichlet boundary condition. We
use a time continuation technique that consists in the introduction of a
fictitious time parameter. We use a POD approach and we present some
theoretical results showing how to choose the optimal dimension of the POD
basis. The results of our computations, related to the first eigenvalue,
confirm the optimal behavior of our approximate solution.
\end{abstract}




\section{Introduction}\label{intro}

The mathematical analysis of FEM approximation to PDE eigenvalue problems is
quite mature and its evolution has reached a deep development in most
application areas.
The a priori analysis of standard Galerkin approximations, after the pioneer
works by Fix and Strang ~\cite{StrangFix} is commonly referred to as the
Babuška-Osborn theory ~\cite{BabuskaOsborn}
In the case of standard Galerkin approximations, it is well understood that a
finite element scheme can be used successfully for the discretization of
compact eigenvalue problems, whenever it works well for the discretization of
the corresponding source problem. This is not the case for other formulations,
such as those arising from mixed methods, where the inf-sup conditions are
neither sufficient nor necessary for the design of a spectrally correct
approximation ~\cite{Boffietal97,Boffietal00,Boffietal13}, ~\cite{Boffi10}.

Reduced order models (ROM) are part of a consolidated technique that is
commonly related to the pioneer works of ~\cite{Almoroth78} , ~\cite{NoorPeters80}, ~\cite{FinkRheinboldt83}, ~\cite{Porshing85}. It underwent then a
successful development in the framework of PDEs from the late 90’s where
different approaches have been explored ~\cite{Peterson89}, ~\cite{ItoRavindran98}, ~\cite{KunischVolkwein01,KunischVolkwein02}, ~\cite{Prudetal02a,Prudetal02b},
~\cite{WillcoxPeraire02}, ~\cite{Veroyetal03}, ~\cite{Benneretal03}, ~\cite{Nguyenetal05}, ~\cite{Huynhetal07}, ~\cite{QuarteroniRozza07}, ~\cite{Rozzaetal08}, ~\cite{Chenetal08}, ~\cite{HaasdonkOhlberger08}, ~\cite{Canutoetal09}, ~\cite{Aguadoetal15}, ~\cite{Hesthavenetal16}.  A commonly used
approach is the Proper Orthogonal Decomposition (POD) method, which is most
often used in the one dimensional temporal discretization as opposed to the
Reduced Basis (RB) approximation which is generally applied to
multi-dimensional parameter domains.

Reduced order model for eigenvalue problems have been the object of a more
limited but not empty investigation, starting from the pioneer works of ~\cite{Madayetal99} and ~\cite{Machielsetal00}, ~\cite{Prudetal02a,Prudetal02b} ,
where only the first fundamental mode was considered. Pau in is PhD thesis
~\cite{Pau07a} has the merit to start the investigation of higher-order modes,
even if still in the setting of simple eigenvalues ~\cite{Pau07b,Pau08}. ~\cite{Buchanetal13} introduced the idea of adding a fictitious time variable for the
solution of a non parametric eigenvalue problem, see also ~\cite{GermanRagusa19}. ~\cite{Fumagallietal16} present an a posteriori analysis for the first
eigenvalue and provide examples of applications to parametric problems.
Finally, ~\cite{Horgeretal17} attack the most complex situation of multiple
eigenvalues by approximating at once a fixed number of eigenmodes, starting
from the first one.

In this paper we investigate the use of reduced order models for the
approximation of eigenvalue problem. We consider the model eigenproblem
associated with the Laplace equation with Dirichlet boundary conditions and we
introduce a fictitious time parameter in order to employ a time continuation
technique. In this framework, we adopt a POD technique in order to build a ROM
after the computation of some snapshots of the high fidelity solution. For the
moment we limit our analysis to the computation of the first fundamental mode.
Future investigations will concern the approximation of other eigenvalues and
of parameter dependent eigenvalue problems.

Section~\ref{se:model} introduces the problem and presents the time
continuation technique approach. Section~\ref{se:rom} describes the reduced
order model framework that we will adopt, the construction of the POD basis,
and the theoretical results about the choice of the optimal dimension for the
POD space. Finally, Section~\ref{se:numres} reports the results of several
numerical tests, showing the good performance of the method.
 
\section{Model problem and time continuation approach}
\label{se:model}

Let us consider the eigenvalue problem related to the Laplace equation with Dirichlet boundary condition: find eigenvalues $\lambda$ and non vanishing eigenfunctions $u$ that solve the following problem

\begin{equation}
\label{orieqn}
\aligned
&-\Delta u(\bfx)=\lambda u(\bfx) &&\text{in }\Omega\\ 
&u(\bfx)=0 &&\text{on }\partial\Omega.
\endaligned
\end{equation}

Our approach considers a reduced order model where a fictitious time parameter $t$ is introduced. In the spirit of ~\cite{Buchanetal13} a time continuation approach, the problem reduces to finding the steady state solution to the following equation 
\begin{equation}\label{teqn}
\frac{\partial u(\bfx,t)}{\partial t}-\Delta u(\bfx,t)=\lambda(t,u) u(\bfx,t).
\end{equation}
This equation must be complemented by a relation linking $u(\bfx,t)$ and $\lambda(t,u)$. In the steady state limit $\lambda$ can be computed by the Rayleigh quotient derived as usual from the stationary equation
\[
-\Delta u(\bfx)=\lambda u(\bfx). 
\]
Multiplying both sides by $u(\bfx)$, integrating over $\Omega$, and using Green's formula  we get
\begin{equation}\label{Rquotient}
\lambda
=\frac{\int_{\Omega} \nabla u(\bfx)\cdot\nabla u(\bfx)\,d\bfx}{\int_{\Omega}u^2 (\bfx)\,d\bfx}=\frac{\int_{\Omega} |\nabla u(\bfx)|^2\,d\bfx}{\int_{\Omega}u^2 (\bfx)\,d\bfx}.
\end{equation}

Hence, the eigenvalues and eigenfunctions of~\eqref{orieqn} are given by the steady state solutions of the following coupled problem: given an initial guess $u^0:\Omega\to\RE$, for all $t>0$ find $u(\bfx,t):\Omega\to\RE$ and $\lambda(t,u)\in\RE$ such that
\begin{equation}
\left\{
\aligned
&\lambda(t,u)=\frac{\int_{\Omega} |\nabla u(\bfx,t)|^2\,d\bfx}{\int_{\Omega}u^2(\bfx,t)\,d\bfx}\\
&\frac{\partial u(\bfx,t)}{\partial t}-\Delta u(\bfx,t)=\lambda(t,u) u(\bfx,t)&&\text{in }\Omega\\
&u(\bfx,0)=u^0(\bfx)&&\text{in }\Omega\\
&u(\bfx,t)=0&&\text{on }\partial\Omega,\ t>0.
\endaligned
\right.
\label{eq:timecont}
\end{equation}

\subsection{Finite element discretization}
We solve Problem~\eqref{eq:timecont} with finite elements in space and finite difference in time. The time derivative is discretized by Euler method, with time step $\Delta t$, as $u_t=\frac{u^{k+1}-u^k}{\Delta t}$, so that we get the following time semi-discretized model
\begin{equation}\label{tdiseqn}
\left\{
\aligned
&\lambda^k=\frac{\int_{\Omega} |\nabla u^k|^2\,d\bfx}{\int_{\Omega}|u^k|^2\,d\bfx}\\
&\frac{1}{\Delta t}u^{k+1}(\bfx)-\Delta u^{k+1}(\bfx)=\lambda^k u^{k}(\bfx)+\frac{1}{\Delta t}u^{k}(\bfx),
\endaligned
\right.
\end{equation}
where $\lambda^k$ is calculated from the Rayleigh quotient \eqref{Rquotient} using the solution $u^k$ at the previous time step.
\begin{figure*}
 \centering
  \includegraphics[width=12cm]{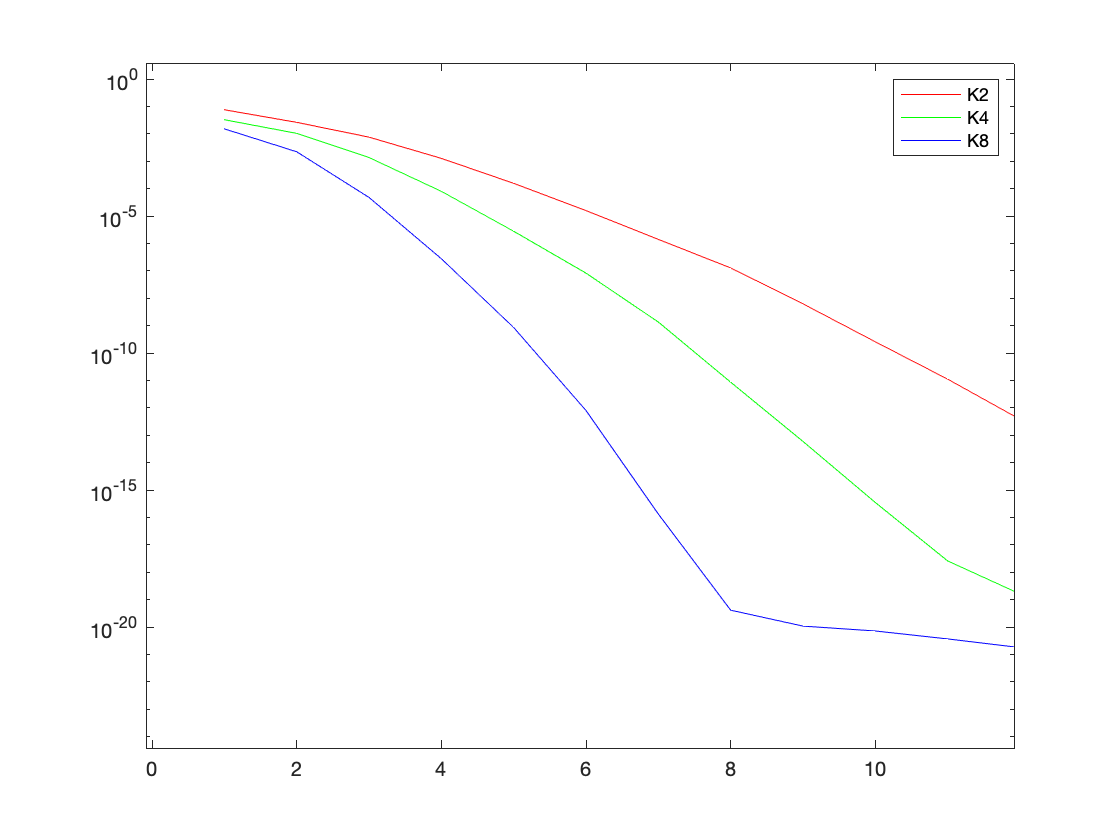}
  \caption{Singular values corresponding to the snapshot matrices $\mathbb{S}_2,\mathbb{S}_4$ and $\mathbb{S}_8$ for crisscross mesh with n=16.}
      \label{fig2}       
 \end{figure*} 

The weak formulation of the first equation in~\eqref{tdiseqn} is obtained as usual after multiplication by a test function $v\in V:=H_0^1(\Omega)$ and integration over $\Omega$
\begin{align*}
\frac{1}{\Delta t} \int_{\Omega}  u^{k+1}v\,d\bfx- \int_{\Omega}  \Delta u^{k+1} v\,d\bfx=\lambda^k  \int_{\Omega}  u^{k}v\,d\bfx+  \frac{1}{\Delta t}\int_{\Omega} u^{k}v\,d\bfx
\end{align*}
Applying Green's theorem in the second term we get
\begin{align*}
\frac{1}{\Delta t} \int_{\Omega}  u^{k+1}v\,d\bfx- \int_{\Omega}  \nabla u^{k+1}\cdot\nabla v\,d\bfx=\lambda^k  \int_{\Omega}  u^{k}v\,d\bfx+  \frac{1}{\Delta t}\int_{\Omega} u^{k}v\,d\bfx
\end{align*}
where we have used the fact that $v=0$ on the boundary. The Galerkin method considers a finite dimensional subspace $V_h \subset V$ of dimension $N_h$ so that the fully discretized problem reads: given an approximation $u^0_h$ of $u^0$, for all $k\ge0$ find $\lambda^k_h\in\RE$ and $u^{k+1}_h\in V_h$ such that
\[
\left\{
\aligned
&\lambda^k_h=\frac{\int_{\Omega} |\nabla u^k_h|^2\,d\bfx}{\int_{\Omega}|u^k_h|^2(\bfx,t)\,d\bfx}\\
&\frac{1}{\Delta t} \int_{\Omega}  u^{k+1}_hv\,d\bfx- \int_{\Omega}  \nabla u^{k+1}_h\cdot\nabla v\,d\bfx=\lambda^k_h  \int_{\Omega}  u^{k}_hv\,d\bfx+  \frac{1}{\Delta t}\int_{\Omega} u^{k}_hv\,d\bfx&&\forall v\in V_h.
\endaligned
\right.
\]

Let $\varphi_i(\bfx)$, $i=1,2,\dots,N_h$, be a basis of $V_h$ and consider the representation $u^k_h(\bfx)=\sum\limits_j u_j^{k}\varphi_j(\bfx)$. Choosing $v=\varphi_i$ for $i=1,2,\dots, N_h$, we finally get the following discrete form of the evolution equation
\[
\aligned
&\frac{1}{\Delta t} \sum_j u_j^{k+1}\int_{\Omega}\varphi_j (\bfx)\varphi_i(\bfx)\,d\bfx+\sum_j u_j^{k+1}\int_{\Omega}\nabla\varphi_j(\bfx)\cdot\nabla\varphi_i(\bfx)\\ &\qquad=\left(\lambda^k_h+\frac{1}{\Delta t}\right)\sum_ju_j^{k}\int_{\Omega}\varphi_j(\bfx)\varphi_i(\bfx).
\endaligned
\]
The matrix form of the evolution problem is then
\begin{equation}\label{matrixform}
AU^{k+1}+\frac{1}{\Delta t} MU^{k+1}=\left(\lambda^k_h+\frac{1}{\Delta t}\right) M U^k
\end{equation}
where
$$
M_{i,j}=\int_{\Omega}\varphi_i(x)\varphi_j(x)dx\quad  \text{and} \quad A_{i,j}=\int_{\Omega}\nabla \varphi_i(x)\cdot \nabla \varphi_j(x)dx.
$$ 
Here $U^k$ denotes the vector of coefficients at the time step $k$. We can rewrite the formula for calculating eigenvalues using these two matrices in discrete form as follows:
\[
\lambda^k_h=\frac{(U^k)^\top AU^k}{(U^k)^\top MU^k}
\]
where $U^\top$ is the transpose of $U$.
Starting from some initial guess $U^0$ and using the equation above we can evaluate the solution until the steady state is reached.
\section{Reduced order model}
\label{se:rom}

With a reduced order model we aim at calculating an approximate solution to the discrete problem \eqref{matrixform} in a subspace $V_N\subset V_h$ of dimension $N$, where $N\ll N_h$. The definition of the subspace $V_N$ is related to the choice of a suitable basis which stems from the high fidelity solution corresponding to an appropriate choice of the parameters. Possible approaches for calculating the basis of $V_N$ include the use of a greedy approach or of a Proper Orthogonal Decomposition (POD). In the greedy approach only $N$ snapshots are selected based on some optimal criterion. On the other hand, in the POD approach $n_s$ snapshots are taken at $n_s$ different predefined values of the parameters. Different sampling techniques like tensorial sampling, Monte Carlo sampling, Latin hypercube sampling, Clenshaw-Curits points may be used.

Here we shall adopt a POD approach.
Let $\mathbb{S}$ be the snapshot matrix of size $N_h\times n_s$, that is the $n_s$ columns of the $\mathbb{S}$ are formed by the snapshots computed with the high fidelity model. The singular values of $\mathbb{S}$ are computed and sorted in decreasing order. Then the first $N$ left singular vectors are taken as the basis of the subspace $V_N$. In the next subsection we are discussing how to calculate the first $N$ left singular vectors of $\mathbb{S}$.

\subsection{Construction of the POD basis}

Let $\mathbb{S}$ be the matrix of the $n_s$ snapshots. The SVD gives
$$\mathbb{S}=\mathbb{U}\Sigma \mathbb{Z}^\top$$
with $$\mathbb{U}=[\pmb{\zeta}_1,\dots,\pmb{\zeta}_{N_h}] \in \mathbb{R}^{N_h\times N_h} \quad \text{and} \quad \mathbb{Z}=[\pmb{\psi}_1,\cdots,\pmb{\psi}_{n_s}] \in \mathbb{R}^{n_s\times n_s}$$
orthogonal matrices and $\Sigma=\mathrm{diag}(\sigma_1,\dots,\sigma_r,0,\dots,0)\in \mathbb{R}^{N_h \times n_s}$ with $\sigma_1\geq \sigma_2\geq \cdots \geq \sigma_r>0$, where $r$ is the rank of the matrix $\mathbb{S}$. As we mentioned earlier, $\{\pmb{\zeta}_1,\dots,\pmb{\zeta}_{N}\}$ will be the basis of the low-dimensional subspace $V_N$. We will not use the SVD directly to find the left singular vectors. Next, we describe the process of calculating the POD basis.
We can write
\begin{align}
\mathbb{S} \pmb{\psi}_i=\sigma_i\pmb{\zeta}_i \quad \text{and}\quad \mathbb{S}^\top \pmb{\zeta}_i =\sigma_i\pmb{\psi}_i, \quad i=1,\dots,r
\end{align}
or, equivalently,
\begin{align}
\mathbb{S}^\top\mathbb{S} \pmb{\psi}_i=\sigma^2_i\pmb{\psi}_i \quad \text{and}\quad \mathbb{S}^\top\mathbb{S} \pmb{\zeta}_i =\sigma^2_i\pmb{\zeta}_i, \quad i=1,\dots,r.
\end{align}
The matrix $\mathbb{C}=\mathbb{S}^T\mathbb{S}$ is called correlation matrix. The POD basis $\mathbb{V}$ of dimension $N\leq n_s$ is defined as the set of the first $N$ left singular vectors $\{\pmb{\zeta}_1,\dots,\pmb{\zeta}_N\}$ of $\mathbb{S}$ or, equivalently, the set of vectors
\begin{align}
\pmb{\zeta}_j=\frac{1}{\sigma_j} \mathbb{S} \pmb{\psi}_j, \qquad 1\leq j\leq N
\end{align}
obtained from the first $N$ eigenvectors $\{\pmb{\psi}_1,\dots,\pmb{\psi}_N\}$ of the correlation matrix $\mathbb{C}$. So, instead of using the SVD of the matrix $\mathbb{S}$ we solve the eigenvalue problem $\mathbb{C} \pmb{\psi}=\sigma^2 \pmb{\psi}$ which is of size $n_s$; notice that $n_s$ is usually much smaller than $N_h$.

Let us now discuss the relation between the solution of the problem in $V_h$ and the approximate solution in $V_N$. Since the basis functions $\pmb{\zeta}_m$, $m=1,2, \dots N $  belong to $V_h$, they can be expressed as linear combinations of basis of $V_h$, namely
$$
\pmb{\zeta}_m=\sum\limits_{i=1}^{N_h}\zeta_m^i \pmb{\varphi}_i \quad 1\leq m\leq N.
$$
The transformation matrix $\mathbb{V}$ is actually the matrix containing all the coefficients $\zeta_m^i$, that is
$$
(\mathbb{V})_{i,m}=\zeta_m^i \quad 1\leq m \leq N, 1\leq i\leq N_h.
$$
Hence, for any linear function $f$, we have $f(\pmb{\zeta}_m)= \sum\limits_{i=1}^{N_h}\zeta_m^i f(\pmb{\varphi}_i)$, $1\leq m\leq N$ or, in matrix form, $\textbf{f}_N=\mathbb{V}^\top\textbf{f}_h$ where $(\textbf{f}_h)_i=f(\pmb{\varphi}_i)$ and $(\textbf{f}_N)_m=f(\pmb{\zeta}_m)$. Since $\mathbb{V}$ is orthogonal we have $\textbf{f}_h=\mathbb{V}\textbf{f}_N$.
The solution $U^k$ of Equation~\eqref{matrixform} is approximated by $U^k=\mathbb{V}{U}_N^k$, where $U_N^k$ is the vector containing the coefficients of the approximate solution at time level $k$ in the subspace $V_N$.
Finally, premultiplying Equation~\eqref{matrixform} by $\mathbb{V}^\top$ and using the relation $U^k=\mathbb{V}{U_N}^k$, we get the formula for reduced order solution at time level $k+1$ as:
\begin{align*}
\mathbb{V}^\top A \mathbb{V}{U}_N^{k+1}+\frac{1}{\Delta t}\mathbb{V}^\top M \mathbb{V}{U}_N^{k}= \left(\lambda^k+ \frac{1}{\Delta t}\right)\mathbb{V}^\top M \mathbb{V}{U}_N^{k}.
\end{align*}
The approximate solution in $V_h$ is given by $U^{k+1}=\mathbb{V}{U}_N^{k+1}$ and the corresponding eigenvalue can be calculated using \eqref{Rquotient}.

\subsection{POD-based RB methods and selection of the dimension of the basis}

The POD basis is orthonormal and it minimizes the sum of squares of errors between each snapshot vector $\pmb{u}_i$ and its projection onto any $N$-dimensional subspace $W$. More precisely it satisfies the following proposition which is a consequence of Schmidt-Eckart-Young theorem.

\begin{theorem}[Schmidt-Eckart-Young] [see~\cite{Quarteronietal16}] \label{SEY}
Given a matrix $\mathbb{A}\in \mathbb{R}^{m\times n}$ of rank r, the matrix $$ \mathbb{A}_k=\sum\limits_{i=1}^k \sigma_i \pmb{\zeta}_{i} \pmb{\psi}_i^T, \quad 1\leq k\leq r $$ satisfies the optimality property
\begin{equation}
\|\mathbb{A}-\mathbb{A}_k\|_F=\min\limits_{\mathbb{B}\in \mathbb{R}^{m\times n}, rank(\mathbb{B})\leq k} \|\mathbb{A}-\mathbb{B}_k \|_F =\sqrt{\sum\limits_{i=k+1}^r \sigma^2},
\end{equation}
where $\|\cdot\|_F$ is the Frobenious matrix norm.
\end{theorem}
\begin{proposition}[see~\cite{Quarteronietal16}]\label{prop1}
Let $\mathscr{V}_N= \{ \mathbb{W}\in \mathbb{R}^{N_h\times N}: \mathbb{W}^T\mathbb{W}=I_N \} $ be the set of all $N$-dimensional bases. Then 
\begin{equation}
\sum\limits_{i=1}^{n_s}\|\mathbf{u}_i-\mathbb{V}^T\mathbb{V} \mathbf{u}_i\|^2=\min\limits_{\mathbb{W}\in \mathscr{V}_N}\sum\limits_{i=1}^{n_s}\|\mathbf{u}_i-\mathbb{W}^T\mathbb{W} \mathbf{u}_i\|^2=\sum\limits_{i=N+1}^{r}\sigma_i^2.
\end{equation}
\end{proposition}

This proposition implies that the error in the POD basis is the sum of the squares of the singular values corresponding to the neglected POD modes. This result gives rise to a criterion for selecting the optimal dimension $N$ of the POD basis. For a given tolerance $\varepsilon$, $N$ will be the smallest integer such that the percentage of energy $I(N)$ of the first $N$ POD modes is less than or equal to $\varepsilon$, that is
\begin{equation}\label{criterion}
I(N)=\frac{\sum\limits_{i=1}^N\sigma_i^2}{\sum\limits_{i=1}^r\sigma_i^2} \geq 1-\varepsilon^2.
\end{equation}

As a typical example of convergence result, we recall the a priori estimate of~\cite{locconv,glbconv}. It refers to the case of elliptic PDEs with single parameter where the associated bilinear form is
\begin{equation}
a(u,v;\mu)=a_0(u,v)+\mu a_1(u,v)\quad \forall u,v\in V
\end{equation}
with $a_0:V \times V\to \mathbb{R}$ and $a_1:V \times V\to \mathbb{R}$ continuous, symmetric, positive semi-definite, and $a_0$ coercive. Then there exists a positive constant $\gamma_1$ such that
$$
0\leq \frac{a_1(v,v)}{a_0(v,v)}\leq \gamma_1 \quad \forall v\in V.
$$
Let us choose the parameter sample points $\mu_i$, $i=1,2,\dots,N$ log-equidistributed in the parameter space $\mathscr{D}$. With this choice of parameters the solution of the reduced order model $u_N(\mu)$ in the space spanned by the corresponding solutions satisfies the following theorem which gives the a priori exponential convergence, where $|||\cdot|||^2=a_0(\cdot,\cdot)$.

\begin{proposition}
For $N\geq N_{crit} $
\begin{align}
||| u_h(\mu)-u_N(\mu)|||\leq (1+\mu_{max}\gamma_1)^{\frac{1}{2}} |||u_h(0)||| e^{-N/N_{crit}}, \quad \forall \mu \in \mathscr{D}
\end{align}
where $N_{crit}=c^{\star} e\cdot log(\gamma \mu_{max}+1)$ with $\gamma$ is an upper bound of $\gamma_1$.
\end{proposition}

In ~\cite{Rozzaetal08} the sampling parameters are chosen as
\begin{equation}
\mu_i=\mu_{min}\exp \Big\{ \frac{i-1}{N-1} \ln(\frac{\mu_{max}}{\mu_{min}})\Big\} \quad 1\leq i\leq N,
\end{equation}
where the parameter space is $\mathscr{D}=[\mu_{min},\mu_{max}]$ and the corresponding low-dimensional space is
$$
W_N=span\{ u_h(\mu_i), 1\leq i\leq N\}.
$$
With this choice, the following results can be obtained.

\begin{proposition}
For any $N\geq N_{crict}$ and for all $\mu \in \mathscr{D}$
$$
\frac{|||u_h(\mu)-u_N(\mu) |||}{|||u_h(\mu) |||} \le\exp \Big\{ -\frac{N-1}{N_{crict}-1}\Big\}
$$
where $N_{crict}=1+[2e\ln\mu_r]$ and $u_h(\mu),u_N(\mu)$ are the FE solution and reduced-order solution in low-dimensional space $W_N$, respectively.
\end{proposition}

\begin{remark}
The error decays exponentially fast with respect to $N$ under the assumption that the solution manifold is analytic~\cite{Quarteronietal16}.
\end{remark}



\section{Numerical results}
\label{se:numres}

 In this section we show the results of our computations with the approach described above. In all our tests we are approximating the first eigenvalue of the model problem. We consider different domains with different types of triangular meshes and use a proper orthogonal decomposition (POD) based reduced order model (ROM).
 We consider the square domain $(0,\pi)^2$ and the L-shaped domain $(-1,1)^2\setminus \{[0,1]\times [0,-1] \}$ with various mesh sequences: uniform meshes of type crisscross, right, and left, and non uniform meshes of Delaunay type.

In our case, the parameter is the time $t$. We consider the time step size $\Delta t=0.1$ and set the stopping criterion as $\frac{\|u^{k+1}-u^{k}\|}{\|u^{k+1}\|} \leq 10^{-8}$ for both ROM and FOM. We use three sets of equally spaced parameters to form snapshot matrices. The corresponding snapshot matrices are denoted by $\mathbb{S}_2, \mathbb{S}_4, \mathbb{S}_8$, respectively, as they consists of solutions of every $2, 4$ and $8$ time steps. We developed our FEM code based on the code provided in ~\cite{LarsonBengzon13,FunkenSchmidt19}.

\begin{figure}
 \centering
 \includegraphics[height=5cm,width=5cm]{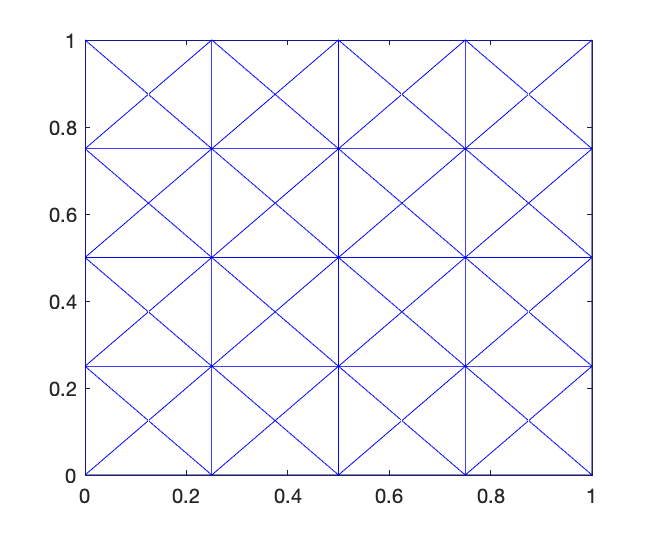}
 \includegraphics[height=5cm,width=5cm]{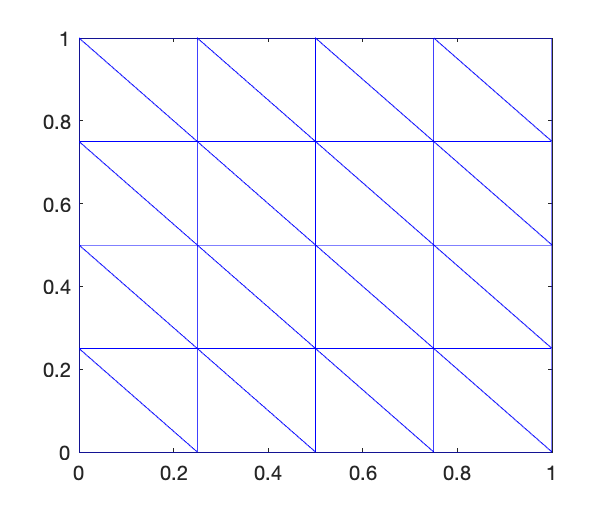}
 \includegraphics[height=5cm,width=5cm]{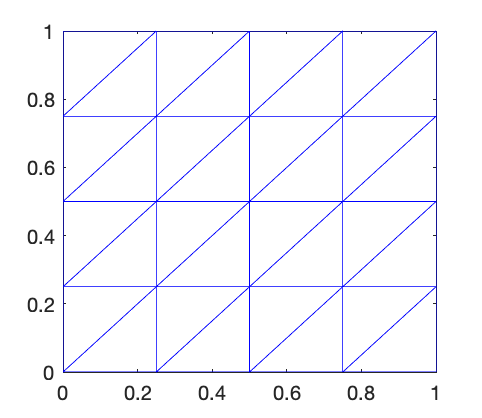}
 \includegraphics[height=5cm,width=5cm]{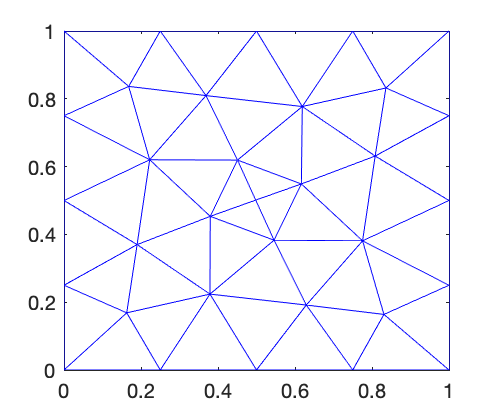}
  \caption{Different type of triangulation on unit square mesh with  n=4 sub-interval on each side. Crisscross, Right, Left and Delaunay triangulation.}
     \label{fig1}       
 \end{figure} 
\subsection{Results on the square domain}

First, we compare the different choices of snapshot matrices mentioned above. Table~\ref{table1}, shows the first eigenvalue corresponding to the different snapshot matrices on the crisscross mesh sequence. We have used P1 FEM for solving the high-fidelity problem. In the table, $n$ denotes the number of sub-intervals on each side of the square domain. The second, third, and fourth columns contain the eigenvalues corresponding to the solution of reduced order model with snapshot matrices $\mathbb{S}_2, \mathbb{S}_4, \mathbb{S}_8$, respectively. In the last column we have reported the minimum number $N$ of POD basis functions satisfying the criterion coming from~\eqref{criterion}, where the tolerance $\varepsilon$ is set as the norm of the difference between the exact and approximate first eigenvectors in normal form. We can see that if the mesh is fine enough then all results corresponding to different snapshot matrices are similar to each other. Also, notice that the results corresponding to the snapshot matrices $\mathbb{S}_4$ and $\mathbb{S}_8$ are similar except the first two cases. From the last column one can see that the results corresponding to $\mathbb{S}_4$ and $\mathbb{S}_8$ are obtained with a smaller number of POD basis functions with respect to $\mathbb{S}_2$ and that
for these two cases the number of POD basis functions is comparable. For this reason, we will use $\mathbb{S}_4$ snapshot matrix in the next results.
In Figure~\ref{fig2} we have shown the first few singular values for the crisscross mesh with number of sub-intervals $n=16$ on each side of the square domain. We can see that the singular values corresponding to $\mathbb{S}_8$ decay much faster than the other two cases. From the theory it is then expected that this snapshot matrix can give better results; this is confirmed by the numbers shown in Table~\ref{table1}.

\begin{table*}[t]
 	 	\centering
 	\begin{tabular}{|c|c|c|c|c|c|c|c|} 
 		\hline
 		n &
		  {\begin{tabular}[c]{@{}c@{}} $\mathbb{S}_2$ \end{tabular}} &
 		 {\begin{tabular}[c]{@{}c@{}}$\mathbb{S}_4$ \end{tabular}} & 
 		{\begin{tabular}[c]{@{}c@{}} $\mathbb{S}_8$  \end{tabular}}&  
 		{\begin{tabular}[c]{@{}c@{}} $N$ \end{tabular}}\\
 	\hline 
     8        & 2.0215691669      &2.0215680613& 2.0215680546&4,4,4\\ 
    16& 2.0053640028(2.1)& 2.0053639952(2.1)&2.0053639952(2.1) &6,5,4\\ 
    32&2.0013392391(2.0)& 2.0013392384(2.0)& 2.0013392384(2.0)&6,5,4\\ 
   64& 2.0003346995(2.0)& 2.0003346995(2.0)&2.0003346994(2.0)& 7,5,5	\\
  128& 2.0000836680(2.0)&2.0000836680(2.0)&2.0000836680(2.0)&8,6,5\\ 
  256&2.0000209166(2.0) & 2.0000209166 (2.0)& 2.0000209166(2.0)   &8,7,5\\ 
  512&2.0000052291(2.0) & 2.0000052291(2.0) & 2.0000052291(2.0)   &9,7,6\\ 
  \hline
 	\end{tabular}
	\caption{ Comparison of first eigenvalue using ROM for three different snapshot matrices on crisscross mesh. For high-fidelity problem we have used P1 space.}
	 	\label{table1}
 \end{table*}
\begin{table*}[t]
 	 	\centering
 	\begin{tabular}{|c|c|c|c|c|c|c|c|} 
 		\hline
 		Mesh &  {\begin{tabular}[c]{@{}c@{}}n \end{tabular}} &
		  {\begin{tabular}[c]{@{}c@{}} FOM \end{tabular}} &
 		 {\begin{tabular}[c]{@{}c@{}}ROM \end{tabular}} & 
		{\begin{tabular}[c]{@{}c@{}} DOF\end{tabular}} & 
 		{\begin{tabular}[c]{@{}c@{}} $N$ \end{tabular}}\\
 	\hline 
Crisscross
&16&2.005363995049&	2.005363995229 &	545	&5	\\ 
&32 &2.001339238351(2.1) &	2.001339238375(2.1) &	2113 &	5\\ 
&64&2.000334699425(2.0) &2.000334699426(2.0)&	8321 &	5\\
&128	&2.000083667969(2.0)&	2.000083667969(2.0)&33025	&6\\
&256	&2.000020916562(2.0)&	2.000020916562(2.0)&131585&	7\\
&512&	2.000005229114(2.0)&	2.000005229114(2.0)&	525313&	7\\

 \hline
Right 
&16&2.019309896556&	2.019309896696	&289	&5\\ 
&32&2.004821215327(2.1)&	2.004821215369(2.1)&	1089	&5\\
&64&2.001204915048(2.1)&	2.001204915048(2.1)&	4225 &	6\\ 
&128&	2.000301204505(2.0)&	2.000301204506(2.0)&	16641&	6\\
&256&	2.000075299611(2.0)&	2.000075299611(2.0)&		66049&	7\\
&512	&2.000018824887(2.0)&	2.000018824807(2.0)&		263169	&7\\
\hline
Left 
&16&2.019309896556&	2.019309896671&	289&	5\\
&32&	2.004821215327(2.1)&	2.004821215368(2.1)&	1089	&5	\\
&64&	2.001204915048(2.1)&	2.001204915048(2.1)&	4225	&6\\ 
&128&	2.000301204505(2.0)&	2.000301204506(2.0)&		16641&	6\\
&256&	2.000075299611(2.0)&	2.000075299611(2.0)&	  66049&	7\\
&512	& 2.000018824807(2.0)&	2.000018824807(2.0)&		263169&	7\\
	\hline
Delaunay 
&16&2.007564672341&	2.007564672491&	467&	5\\ 
&32&	2.001871138383(2.1)&	2.001871138450(2.1)&	1784	 &5\\
&64&	2.000453161338(2.0)&	2.000453161343(2.0)&	7140 &	5\\
&128	&  2.000113145511(2.0)&	2.000113145511(2.0)&		28523&	6\\
&256	& 2.000028027418(2.0)&	2.000028027418(2.0)&	113846&	7\\
&512    & 2.000006988690(2.0)         &2.000006988690(2.0)& 454733 & 7\\
\hline
\end{tabular}
	\caption{ Comparison of first eigenvalue on different meshes using P1 FEM.}
 	\label{table2}
 \end{table*}
 
  \begin{table*}[t]
 	 	\centering
 	\begin{tabular}{|c|c|c|c|c|c|c|c|c|} 
 		\hline
 	Mesh &
		 {\begin{tabular}[c]{@{}c@{}} Method \end{tabular}} &
		  {\begin{tabular}[c]{@{}c@{}} n=16 \end{tabular}} &
 		 {\begin{tabular}[c]{@{}c@{}}n=32 \end{tabular}} & 
 		{\begin{tabular}[c]{@{}c@{}} n=64  \end{tabular}}&  
		{\begin{tabular}[c]{@{}c@{}} n=128\end{tabular}} &
		{\begin{tabular}[c]{@{}c@{}} n=256\end{tabular}} &
 		{\begin{tabular}[c]{@{}c@{}} $512$ \end{tabular}}\\
 	\hline 
	Crisscross & FOM&0.04&0.14&0.85&3.57&16.38& 83.61\\
       & ROM&0.008&0.01&0.03&0.13&0.60&2.79\\
	    \hline
	    Right& FOM&0.02&0.08&0.39&1.66&7.56&42.53\\ 
	             &ROM&0.006&0.008&0.02&0.06&0.24&1.25\\ 
	\hline    
	Left&FOM&0.02& 0.08& 0.39&1.58&7.26 &41.15\\ 
	         &ROM&0.006&0.008&0.02& 0.06&0.25&1.23\\ 
	     \hline
	      Delaunay&FOM & 0.04&0.18&0.90&4.10&22.34& 185.56\\
	                     &ROM& 0.005&0.01&0.04&0.15&0.76&5.52\\
	  \hline              
  \end{tabular}
	\caption{ CPU time for FOM and ROM corresponding to the results of P1 FEM on square domain.}
 	\label{table22}
 \end{table*}
 
 Another interesting subject of investigation is related to the dependence of the results one the different types of triangular meshes. Here we consider our four types of triangular meshes: crisscross, right, left, and Delaunay mesh. A sample mesh on unit square domain has been shown in Figure~\ref{fig1}. As we have already mentioned, we take $\mathbb{S}_4$ as the snapshot matrix also for this test. In Table~\ref{table2} we have reported the first eigenvalue on four different kinds of triangular meshes with different degrees of freedom. In the last column of the table we have indicated the number of POD basis functions which is obtained by utilising the criterion coming from \eqref{criterion}. The rate of convergence of the eigenvalues is written within parentheses. Observe that the rate of convergence is 2 for all the meshes as expected because we are using P1 finite elements and the eigenspace is smooth. The CPU time taken by the full order model and the CPU time for the reduced order model are reported in Table~\ref{table22}. From the CPU time reported in Table~\ref{table22} one can see that in the ROM the online phase is from 10 to 30 times faster than the full order model. Notice that the results on the right and left uniform meshes are exactly the same. These two types of meshes are the same with a different orientation. Hence, from now on we will use only the right uniform mesh for our tests. In all cases, when the mesh is fine enough, then the results obtained by the ROM are exactly same as the ones obtained by the FOM.
 
\begin{table*}[t]
 	 	\centering
 	\begin{tabular}{|c|c|c|c|c|c|c|c|} 
 		\hline
 		Mesh &  {\begin{tabular}[c]{@{}c@{}}n \end{tabular}} &
		  {\begin{tabular}[c]{@{}c@{}} FOM \end{tabular}} &
 		 {\begin{tabular}[c]{@{}c@{}}ROM \end{tabular}} & 
		{\begin{tabular}[c]{@{}c@{}} DOF\end{tabular}} & 
 		{\begin{tabular}[c]{@{}c@{}} $N$ \end{tabular}}\\
 	\hline 
	Crisscross 
	&16&2.0000034983508(4.0) & 2.0000034983774(4.1)&2113&	6\\
	&32&2.0000002191655(4.0)	&2.0000002191679(4.0)	&	8321	&6\\
	&64&2.0000000137061(4.0)&2.0000000137061(4.0)	&33025	&7\\
	&128&2.0000000008569(4.0)&	2.0000000008569(4.0)&	131585	&8\\
	&256&2.0000000000535(4.0)&	2.0000000000535(4.0)&	525313&	8\\
	&512&2.0000000000007(6.0)&	2.0000000000006(6.0)& 2099201	&8\\
	\hline
	Right
	&16&2.0000286902960(4.0)	 &2.0000286902960(4.1)&	1089&6\\
	&32&2.0000018029662(4.0)&2.0000018029662(4.1)&	4225	&7\\
	&64&2.0000001128424(4.0)&	2.0000001128424(4.0)&16641&7\\
	&128&2.0000000070555(4.0)&2.0000000070554(4.0)&66049&7\\
	&256&2.0000000004417(4.0)&2.0000000004417(4.0)&	263169&	7\\
	&512&2.0000000000333(3.7)&	2.0000000000333(3.7)&1050625&8\\
	\hline
	Delaunay
	&16&2.0000058772706(4.1)&	2.0000058772719(4.1)&		1799	&6\\
	&32&2.0000003559407(4.1)&	2.0000003559447(4.1)&	47003	&6\\
	&64&2.0000000216176(4.0)&	2.0000000216177(4.0)&	28299&	7\\
	&128&2.0000000013122(4.0)&	2.0000000013121(4.0)&	113575&7\\
	&256&2.0000000000813(4.0)&	2.0000000000813(4.0)&	454355&8\\
	&512&2.0000000000051(4.0)& 2.0000000000051(4.0)&1816879 &8\\
	\hline
\end{tabular}
	\caption{ Comparison of first eigenvalue on three different mesh with P2 FEM on Square Domain $(0,\pi)^2$.}
 	\label{table3}
 \end{table*}

 \begin{table*}
 \begin{tabular}{|c|c|c|c|c|c|c|c|c|} 
 		\hline
 		Mesh &
		 {\begin{tabular}[c]{@{}c@{}} Method \end{tabular}} &
		  {\begin{tabular}[c]{@{}c@{}} n=16 \end{tabular}} &
		  {\begin{tabular}[c]{@{}c@{}} n=32 \end{tabular}} &
 		 {\begin{tabular}[c]{@{}c@{}}n=64 \end{tabular}} & 
 		{\begin{tabular}[c]{@{}c@{}} n=128  \end{tabular}}&  
		{\begin{tabular}[c]{@{}c@{}} n=256\end{tabular}} &
		{\begin{tabular}[c]{@{}c@{}} n=512\end{tabular}} \\
 	\hline 
	Crisscross & FOM&0.16 &0.78&3.40&17.53&86.37&414.58\\
	                 & ROM&0.02&0.04&0.17  &0.89&3.80&22.12\\
	    \hline
	    Right& FOM &0.10& 0.47&1.80 &8.80&46.65&222.04\\
	             &ROM&0.01&0.02&0.08&0.39&1.74&15.52  \\
	                \hline
	      Delaunay&FOM &0.20&0.05 &4.60&24.30&145.71& 847.67\\
	                     &ROM&0.02&0.04&0.19&1.29&5.12&38.40\\
	  \hline              
  \end{tabular}
	\caption{ CPU time for FOM and ROM corresponding to the results of P2-FEM on square domain.}
 	\label{table31}
 \end{table*} 

 \begin{figure*}
 \centering
 \includegraphics[height=4cm,width=4cm]{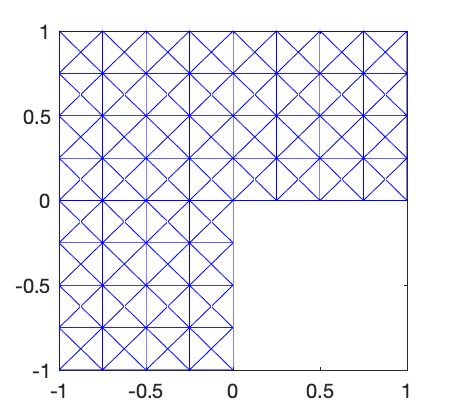}
 \includegraphics[height=4cm,width=4cm]{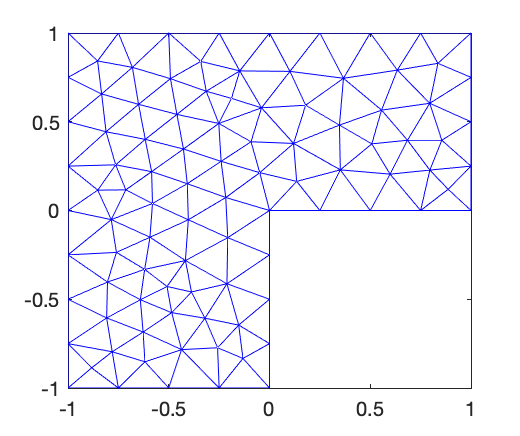}
 \includegraphics[height=4cm,width=3.9cm]{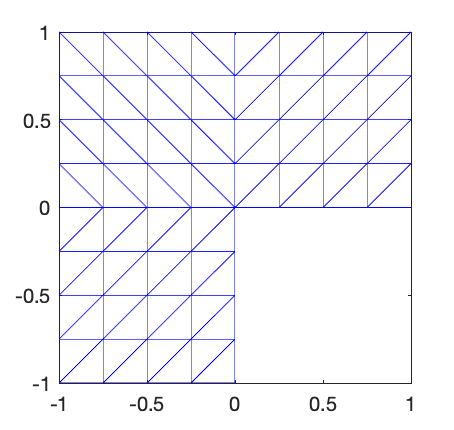}
  \caption{Different type of Triangular mesh with n=4 subinterval on each side of the squares. Left to right: Crisscross, Delaunay, Mixed triangulation.}
     \label{fig3}     
\end{figure*}  
       
Finally, we have used P2 FEM for the high-fidelity problem in the offline stage. The results corresponding to P2 FEM on the square meshes are reported in Table~\ref{table3}. In this case we have considered triangular meshes of type crisscross, right, and Delaunay. For crisscross mesh we took a number of sub-intervals equal to $8$, $16,\dots$, $512$, that is we have done uniform refinement in each step. The results are converging to the exact solution with rate of convergence equal to $4$. In the last case the result is correct up to six decimal place. Like in the P1 case, here also we have reported the CPU time for both the FOM and the ROM computations in Table~\ref{table31}. From the table it is clear that the ROM is 10-15 times faster than the corresponding FOM. Similarly, for the right uniform mesh also the order of convergence is equal to $4$ except for the last case. In this case also the ROM is giving the solution 10-20 times faster than the offline stage. In the finest mesh the eigenvalue is matching the exact one up to 10 decimal place and the corresponding number of degrees of freedom is $1,050,625$. Similar conclusions can be drawn for the unstructured Delaunay mesh.

 \begin{table*}[t]
 	 	\centering
 	\begin{tabular}{|c|c|c|c|c|c|c|c|} 
 		\hline
 		Mesh &  {\begin{tabular}[c]{@{}c@{}}n \end{tabular}} &
		  {\begin{tabular}[c]{@{}c@{}} FOM \end{tabular}} &
 		 {\begin{tabular}[c]{@{}c@{}}ROM \end{tabular}} & 
 		{\begin{tabular}[c]{@{}c@{}} DOF\end{tabular}} & 
 		{\begin{tabular}[c]{@{}c@{}} $N$ \end{tabular}}\\
  	\hline 
Crisscross 
&16&	9.689550909492&	9.689550909499&	1601	&5\\ 
cross&32&    9.657193968400(1.53)&	9.657193968400(1.53)&	6273	& 6\\
&64&	9.646089739637(1.47)&	9.646089739637(1.47)&	24833&	6\\
&128	&    9.642109199667(1.42)&	 9.642109199667(1.42)&		98817&	6\\
&256&	9.640635313503(1.39)&	9.640635313503(1.39)&	394241	&7	\\
&512	 & 9.640076778359(1.37)&	9.640076778359(1.37)&		1574913&	7\\
 \hline
 Mixed
 &16&	9.737622827078&	9.737622827084&		833&	5\\
 &32&	9.672057256699(1.65)&	9.672057256699(1.65)&		3201	&6\\
 &64& 	9.650967866448(1.55)&	9.650967866448(1.55)&	12545&	6\\
 &128&	9.643793757105(1.48)&	9.643793757105(1.48)&		49665&	6\\
 &256&	9.641241417249(1.43)&	9.641241417249(1.43)&	197633&	7\\
 &512&	9.640301753187(1.40)&	9.640301753187(1.40)&	788481&	7	\\
 \hline
 Delaunay
 &16& 	9.691713765278&	9.691713765278&		1356	& 5\\
 &32&	9.659203788910(1.43)&	9.659203788910(1.43)&		5354	 &6	\\
 &64&	9.645755044694(1.69)&	9.645755044694(1.69)&	21406	&6	\\
 &128&	9.641949706951(1.44)&	9.641949706951(1.44)&	85651&	6\\
 &256&	9.640665163096(1.25)&	9.640665163096(1.25)&		338736	&7\\
 
 \hline
\end{tabular}
	\caption{ Comparison of first eigenvalue computed by ROM using three different mesh with P1-FEM On L-shape Domain.}
 	\label{table4}
 \end{table*}
 
 \begin{table*}
 \begin{tabular}{|c|c|c|c|c|c|c|c|c|} 
 		\hline
 		Mesh &
		 {\begin{tabular}[c]{@{}c@{}} Method \end{tabular}} &
		  {\begin{tabular}[c]{@{}c@{}} n=16 \end{tabular}} &
		  {\begin{tabular}[c]{@{}c@{}} n=32 \end{tabular}} &
 		 {\begin{tabular}[c]{@{}c@{}}n=64 \end{tabular}} & 
 		{\begin{tabular}[c]{@{}c@{}} n=128  \end{tabular}}&  
		{\begin{tabular}[c]{@{}c@{}} n=256\end{tabular}} &
		{\begin{tabular}[c]{@{}c@{}} n=512\end{tabular}} \\
 	\hline 
	Crisscross & FOM&0.9&0.40&2.07&9.69&48.94&226.98\\
	          & ROM&0.01&0.03&0.10&0.45&2.31&9.68\\
	 	  \hline    
	Right&FOM&0.05& 0.22&1.09&5.15&23.60&124.88\\
	        &ROM&0.008&	0.01& 0.04&0.22&1.07&4.60\\
		  \hline   
Delaunay&FOM&0.10&0.59&2.63&13.47&109.02&-\\
               & ROM& 0.01&0.03&0.11&0.57&3.32&- \\
		  \hline       
  \end{tabular}
	\caption{ CPU time for FOM and ROM corresponding to the results of P1-FEM on L-shape domain.}
 	\label{table41}
 \end{table*}

\subsection{Results on the L-shaped domain}

In this last section we investigate how our ROM method is performing on the L-shaped domain $(-1,1)^2\setminus \{[0,1]\times [0,-1] \}$ where singular solutions are expected. We consider three types of triangular meshes: crisscross, Delaunay and what we call mixed uniform mesh, which is a combination of the right and left meshes. A sample mesh has been shown in Figure~\ref{fig3}. In this case $n=4$ is the number of sub-intervals in each side of the unit square and the L-shaped domain consists of three unit squares.

In Table~\ref{table4} we have reported the results of our reduced order model on the singular L-shaped domain using three different triangular meshes and P1 FEM in the full order model. The rate of convergence is reported within parentheses. We have reported the CPU time for FOM and ROM in Table~\ref{table41}. We can see that ROM is at least 13 times faster than FOM. For example, on the crisscross uniform mesh with $n=256$, the time required for FOM is $48.94$ seconds, whereas CPU time for ROM is $2.31$ seconds, that is $21$ times less. The rate of convergence is around $1.4$ where we take as a reference value for the first eigenvalue $9.6397238440219$. It can be seen that the eigenvalue matches up to only $2$ decimal places the exact value; this is as usually due to the singularity present in the domain. Similarly, we can see for the right mesh the order of convergence is approximately $1.4$, but for Delaunay mesh it is varying because it is an unstructured mesh. In all cases we can see that the reduced model is at least $10$ times faster.

\begin{table*}[t]
 	 	\centering
 	\begin{tabular}{|c|c|c|c|c|c|c|c|} 
 		\hline
 		Mesh &  {\begin{tabular}[c]{@{}c@{}}n \end{tabular}} &
		  {\begin{tabular}[c]{@{}c@{}} FOM \end{tabular}} &
	 		 {\begin{tabular}[c]{@{}c@{}}ROM \end{tabular}} & 
 		{\begin{tabular}[c]{@{}c@{}} DOF\end{tabular}} & 
 		{\begin{tabular}[c]{@{}c@{}} $N$ \end{tabular}}\\
 	\hline 
Crisscross
&8& 9.654566452569&	9.654566452611&	1601&	5\\
&16&9.645648130148(1.33)&	9.645648130148(1.33)&		6273	&6\\
&32&9.642079379771(1.33)&	9.642079379771(1.33)&		24833&	6\\
&64&9.640659153042(1.33)&9.640659153042(1.33)&		98817&	6\\
&128&9.640095077363(1.33)&	9.640095077363(1.33)&	394241&	7\\
&256&9.639871173962(1.33)&	9.639871173962(1.33)&	1574913&	7\\
\hline
Mixed
&8& 9.663207242395&	9.663207242479&		833&	5\\
&16&9.649100751840(1.32)&	9.649100751840(1.32)&	3201 &	6\\
&32&9.643459074163(1.33)&	9.643459074163(1.33)&	12545	&6\\
&64&9.641208083352(1.33)&	9.641208083352(1.33)&	49665&	6\\
&128&9.640313088235(1.33)&9.640313088235(1.33)&	197633	&6\\
&256&9.639957710191(1.33)&	9.639957710191(1.33)&	788481&	7\\
\hline
Delaunay
&8& 9.655492013317&9.655492013320	&	1385	 &5\\
&16&9.646108037396(1.30)&	9.646108037396(1.30)&	5293&	6\\
&32&9.642596253224(1.15)&	9.642596253225(1.15)&	21157	&6\\
&64&9.640617826930(1.68)&9.64061782693301.68)&	85109&	6\\
&128&9.640070529727(1.37)&	9.640070529727(1.37)&	341577&	6\\
&256&	9.639886064681(1.10)&	9.639886064681(1.10)& 1352893&	7\\
 \hline
\end{tabular}
	\caption{ Comparison of First Eigenvalue on ROM three different mesh with P2 FEM and Rayleigh Quotient Formula and $\mathbb{S}_4$ as snapshot matrix On L-shape Domain.}
 	\label{table5}
 \end{table*}
 
\begin{table*}
 \begin{tabular}{|c|c|c|c|c|c|c|c|c|} 
 		\hline
 		Mesh &
		 {\begin{tabular}[c]{@{}c@{}} Method \end{tabular}} &
		 {\begin{tabular}[c]{@{}c@{}} n=8\end{tabular}} &
		  {\begin{tabular}[c]{@{}c@{}} n=16 \end{tabular}} &
		  {\begin{tabular}[c]{@{}c@{}} n=32 \end{tabular}} &
 		 {\begin{tabular}[c]{@{}c@{}}n=64 \end{tabular}} & 
 		{\begin{tabular}[c]{@{}c@{}} n=128  \end{tabular}}&  
		{\begin{tabular}[c]{@{}c@{}} n=256\end{tabular}} \\
 	\hline 
	Crisscross & FOM&0.13&0.45 &2.13&10.59&52.57&253.54\\
	                 & ROM& 0.01&0.03&0.14&0.70&3.30 &17.69\\
	    \hline
	    Mixed& FOM &0.077&0.25 &1.34&5.78&27.18&141.82\\
	              &ROM& 0.007&0.02&0.07&0.32&1.47 &6.55\\
	                \hline
	      Delaunay&FOM &0.13&0.74&2.92&15.15&95.08&501.30\\
	                     &ROM &0.01&0.03 &0.15&0.69&3.60&24.35\\
	  \hline              
  \end{tabular}
	\caption{ CPU time for FOM and ROM corresponding to the results of P2-FEM on L-shape domain.}
 	\label{table51}
 \end{table*}
In Table~\ref{table5} we have reported the results of P2 FEM on the L-shaped domain. Here also, we have calculated results on the three types of mesh: crisscross, mixed, and Delaunay. In this case the order of convergence is 1.33 for the crisscross and mixed meshes and varying for the Delaunay mesh because it is a non-structured mesh. In this case the ROM is at least $13$ times faster than the FOM, which can be observed from Table~\ref{table51}, where we have included the CPU time for FOM and ROM, respectively. In the case of the P2 FEM the eigenvalue is approximated correctly up to 3 decimal places.

\begin{table*}[t]
 	 	\centering
 	\begin{tabular}{|c|c|c|c|c|c|c|c|} 
 		\hline
 		Mesh & 
		  {\begin{tabular}[c]{@{}c@{}} FOM \end{tabular}} &
		   {\begin{tabular}[c]{@{}c@{}}ROM \end{tabular}} & 
		{\begin{tabular}[c]{@{}c@{}} DOF\end{tabular}} & 
 		{\begin{tabular}[c]{@{}c@{}} $N$ \end{tabular}}\\
 	\hline 
Crisscross
&9.6397259448(1.64)&	9.6397259448(1.64)&	128876&9\\
&9.6397249480(2.15)&	9.6397249480(2.15)&	173946	&9\\
&9.6397244634(1.76)&	9.6397244634(1.76)&	241486&	10\\
&9.6397241783(1.57)&	9.6397241783(1.57)&	357343&	10\\
&9.6397240778(1.14)&9.6397240778(1.41)&	488670&	11\\
&9.6397239925(1.57)&9.6397239925(1.57)&	652305&	11\\
\hline
Mixed

&9.6397258141(1.92)&9.6397258141(1.92)&112184	&9\\
&9.6397248627(1.63)&9.6397248627(1.63)&168046	&10\\
&9.6397245461(1.16)&9.6397245461(1.16)&231542	&10\\
&9.6397243056(1.44)&	9.6397243056(1.44)&309858&	11\\
&9.6397240956(1.89)&9.6397240956(1.89)&427413&	10\\
&9.6397239677(1.83)&9.6397239677(1.83)&629926&	10\\
\hline
Delaunay
&9.6397309233(1.84)&9.6397309233(1.84)&73219&	8\\
&9.6397272642(1.81)&	9.6397272642(1.81)&109465	&8\\
&9.6397255112	(1.83)&	9.6397255112(1.83)&162254	&8\\
&9.6397247163(1.71)&9.6397247163(1.71)&236753&	9\\
&9.6397243381(1.53)&9.6397243381(1.53)&343460&10\\
&9.6397241186	(1.65)&	9.6397241186(1.65)&	490769&	11\\
\hline
\end{tabular}
	\caption{ Comparison of first eigenvalue on ROM three different meshes with Adaptive P2 FEM on the L-shaped domain.}
 	\label{table6}
 \end{table*}
\begin{table*}
 \begin{tabular}{|c|c|c|c|c|c|c|c|c|} 
 		\hline
 		Mesh &
		 {\begin{tabular}[c]{@{}c@{}} Method \end{tabular}} &
		 {\begin{tabular}[c]{@{}c@{}} I\end{tabular}} &
		  {\begin{tabular}[c]{@{}c@{}} II \end{tabular}} &
		  {\begin{tabular}[c]{@{}c@{}} III \end{tabular}} &
 		 {\begin{tabular}[c]{@{}c@{}}IV \end{tabular}} & 
 		{\begin{tabular}[c]{@{}c@{}} V  \end{tabular}}&  
		{\begin{tabular}[c]{@{}c@{}} VI\end{tabular}} \\
 	\hline 
	Crisscross & FOM &31.84&46.12&74.80&115.27&163.06&230.88\\
	                 & ROM &1.56&1.82&3.16 &5.31&8.10&11.54\\
	    \hline
	    Mixed& FOM &26.25&42.02&67.98&98.71&137.83&218.03\\
	              &ROM &1.08&2.07 &2.90&4.65&6.54&10.80\\
	                \hline
	      Delaunay&FOM &17.08&26.89&44.05 &71.36&121.29&177.53\\
	                     &ROM &0.70 &1.16&1.75&2.79&5.55&8.85 \\
	  \hline              
  \end{tabular}
	\caption{ CPU time for FOM and ROM corresponding to the results of Adaptive P2 FEM on the L-shaped domain.}
 	\label{table61}
 \end{table*}

In order to improve the approximation of the singular eigenfunction, we apply an adaptive scheme based on P2 FEM in order to calculate the solution of the full order model. Using the selected snapshots from the full order model we calculate the approximation of eigenvalues and eigenvectors with the reduced order model. 
For the adaptive method we adopt the usual strategy: Solve, Estimate, Mark, and Refine. The following residual based error estimator is used locally in each element $K$
\begin{equation}\label{residual}
\eta_K^2 =h_K^2\|r|_K\|_{0,K}^2+ \frac{1}{2} \sum\limits_{e\in \varepsilon_K} h_e\|R|_e\|_{0,e}^2
\end{equation} 
with
$$
r|_K=\Delta u_h+\lambda_h u_h \quad \forall K \quad \text{and}\quad R|_e=\nabla u_h |_{K^1_e} \cdot \pmb{n}_{K_e^1}+\nabla u_h|_{K^2_e} \cdot\pmb{n}_{K_e^2}
$$
where $K_e^1$ and $K_e^2$ are the elements sharing the edge $e$ and  $\pmb{n}_{K_e^1}$  and $\pmb{n}_{K_e^1}$ are the outward unit normal vectors on $\partial K_e^1$ and $\partial K_e^1$, respectively.

The results on the crisscross, Delaunay, and mixed meshes are shown in Table~\ref{table6}. We have reported the eigenvalues computed with the FOM in the second column, while the ROM eigenvalues are presented in the third column. In Table~\ref{table61} we have reported the CPU time taken by FOM and ROM for all the results.
In the final refinement we get the approximate eigenvalues equal to $9.6397239925$ for both FOM and ROM on the crisscross mesh, which is quite a good result. In this case we get order of convergence equal to approximately $1.57$. Also the ROM is much much faster than the full order model because in FOM we need to calculate error, refine the mesh and solve the problem. But in ROM we have the mesh stored in FOM and we only solve the reduced problem. 

Similarly, in the final refinement we get an approximate eigenvalues of about $9.6397241186$ on the Delaunay mesh, which is correct up to $5$ decimal digits. In this case the minimum order of convergence is $1.4$. On the mixed uniform mesh the computed eigenvalue is $9.6397239677$, which is correct up to $6$ decimal digits and the minimum order of convergence is $1.4$. In these two cases also the ROM is much much faster than the full order model. In all the cases the number of reduced basis function used is nine or ten.

\section*{Acknowledgements}
This research was supported by the Competitive Research Grants Program CRG2020 ``Synthetic data-driven model reduction methods for modal analysis'' awarded by the King Abdullah University of Science and Technology (KAUST).
Daniele Boffi is member of the INdAM Research group GNCS and his research is
partially supported by IMATI/CNR and by PRIN/MIUR.

\bibliographystyle{plain}
\bibliography{mybib1}

\end{document}